\documentclass[letterpaper, 10 pt, conference]{ieeeconf}  

\IEEEoverridecommandlockouts                              

\usepackage{amsmath,amssymb,amsthm,mathtools}
\usepackage[pdfstartview=FitB,plainpages=false,colorlinks=true,linkcolor=black,citecolor=black,linktocpage=true,hyperfootnotes=false,bookmarksopen=true,urlcolor=black]{hyperref}

\newcommand{\tint}{\textstyle\int}
\newcommand{\Rset}{\mathbb R}
\newcommand{\pprime}{{\prime\prime}}
\newcommand{\dd}{\mathrm{d}}

\newtheorem{remark}{Remark}
\newtheorem{assumption}{Assumption}
\newtheorem{definition}{Definition}

\title{\LARGE \bf
Successive design of backstepping observers for parabolic PDE-ODE systems and its duality to state feedback stabilization
}

\author{Nicole Gehring%
\thanks{Nicole Gehring is with the Chair of Systems Theory and Control Engineering, Otto von Guericke University Magdeburg, Magdeburg, Germany (e-mail: nicole.gehring@ovgu.de).}%
}

\begin{document}

\maketitle
\thispagestyle{empty}
\pagestyle{empty}

\begin{abstract}
The paper introduces a successive backstepping observer design for strictly feedforward parabolic PDE-ODE systems, in which the coupling structure determines the order of error stabilization and the corresponding transformations.
First, a transformation based on a virtual measurement stabilizes the ODE observer error subsystem, which is most distal from the measurement, while decoupling it from the PDE error state.
Second, a Volterra integral transformation is employed to stabilize the PDE error subsystem and to map the overall error dynamics into a cascade of exponentially stable ODE and PDE subsystems.
The design is shown to be dual to a recently proposed multi-step state feedback design for parabolic PDE-ODE systems in strict feedback form, thus explaining the structure of the presented observer design.
\end{abstract}

\section{Introduction}

Backstepping has emerged over the last two decades as a powerful constructive method for control and state estimation of systems governed by partial differential equations (PDEs), with numerous applications to parabolic, hyperbolic, and interconnected systems (see \cite{Krstic2008book} and the survey paper \cite{Vazquez2026aut}).
The first backstepping controller in the seminal paper \cite{Balogh2002ejc} is derived for a boundary-actuated reaction-diffusion equation by considering a spatial discretization of the PDE and building on the recursive backstepping design well-known for ordinary differential equations (ODEs) in strict feedback form (e.g., \cite{Krstic1995}).
The Volterra integral transformation arising from that design is inherently linked to the strict feedback form of the PDE itself.
Based on that and inspired by \cite{Krener2003siam}, \cite{Smyshlyaev2005scl} provides the first backstepping observers for a parabolic PDE with boundary measurements by making use of the duality of controller and observer design for linear systems.

Similar to \cite{Smyshlyaev2005scl} and the generalizations in \cite{Orlov2017siam,Camacho2020scl}, many backstepping observers rely on the strict feedforward form of the system w.r.t.\ the boundary measurement.
For example, this is true for the cascades of ODE and parabolic PDE subsystems studied in \cite[Chapter 17]{Krstic2008book} (see also, e.g., \cite{Vazquez2011CDC,Aamo2013tac,Auriol2020aut} in the hyperbolic case).
PDE-ODE systems with bidirectionally coupled PDE and ODE subsystems are significantly more challenging than the cascaded structures.
They are first addressed in \cite{Tang2011franklin} for a parabolic system with a measurement of the PDE state at the boundary where the two subsystems are interconnected. Such a system is not in strict feedforward form, which is why the design imposes an additional restriction on the location of the eigenvalues of the ODE.
A similar assumption is considered for the heat equation sandwiched between two ODEs in \cite{Wang2019TAC}.

The aforementioned backstepping observers are typically derived using a single transformation, directly constructing an integral transformation that enforces a desired target error dynamics.
For bidirectionally coupled PDE-ODE systems, however, the particular structure of the coupling becomes crucial.
In many physically motivated models, the overall dynamics exhibits a strict feedforward form w.r.t.\ the boundary measurement, which suggests a successive rather than a monolithic observer design.
In fact, the dual structure, the strict feedback form w.r.t.\ a boundary input, has already been exploited for the state feedback design (see again the survey paper \cite{Vazquez2026aut} as well as \cite{Gehring2026ECC}), particularly for multi-step design approaches in \cite{Gehring2021MTNS} and \cite{Deutscher2021TAC}.
Apart from \cite{Deutscher2021TAC}, there are only few observer designs that adopt a multi-step structure.

This paper proposes a successive two-step backstepping observer design for a class of strictly feedforward parabolic PDE-ODE systems.
In contrast to a direct construction of a single observer transformation, the proposed procedure follows the coupling structure of the system.
The ODE error subsystem, which is most distal from the boundary measurement, is stabilized first by introducing a virtual measurement and a dedicated state transformation, while the remaining PDE error dynamics is subsequently stabilized based on a Volterra integral transformation.
The observer transformations and injection functions are shown to be dual to the corresponding transformations and feedback laws of the multi-step state feedback design in \cite{Gehring2021MTNS}.
This explains the structure of the proposed observer approach.
In the end, the design strategy provides a framework that can be extended beyond the specific class of systems considered here to other interconnected systems in strict feedforward form.

Following the problem statement in Section~\ref{sec:problem}, the two-step observer design is presented in Section~\ref{sec:obs}. 
To further explain the design, it is shown in Section~\ref{sec:dual} that this observer design is dual to the state feedback design in \cite{Gehring2021MTNS}.

\section{Problem statement}
\label{sec:problem}

Consider the parabolic PDE-ODE system
\begin{subequations}
	\label{eq:sys}
	\begin{align}
		\label{eq:sys-ode}
		\dot\xi(t) &= F\xi(t) + Bx(0,t) \\
		\label{eq:sys-rb0}
		\partial_z x(0,t) &= Q_0 x(0,t) + C\xi(t) \\
		\partial_t x(z,t) &= \Lambda(z)\partial_z^2 x(z,t) + A(z)x(z,t) \\
		\label{eq:sys-rb1}
		\partial_z x(1,t) &= Q_1 x(1,t) \\
		\label{eq:sys-meas}
		y(t) &= x(1,t)
	\end{align}
\end{subequations}
without an input, where the ODE subsystem \eqref{eq:sys-ode} with state $\xi(t)\in\Rset^{n_0}$ is bidirectionally interconnected with the PDE subsystem \eqref{eq:sys-rb0}--\eqref{eq:sys-rb1} with state $x(z,t)\in\Rset^n$.
The diffusion coefficients $\lambda_i\in C^2([0,1])$, $i=1,\dots,n$, in the diagonal matrix $\Lambda=\mathrm{diag}\,(\lambda_1,\dots,\lambda_n)$ satisfy $\lambda_1(z)>\cdots>\lambda_n(z)>0$, $z\in[0,1]$, the elements of $A(z)\in\Rset^{n\times n}$ are $C([0,1])$ functions.
For simplicity $Q_0\in\Rset^{n\times n}$ is assumed to be a diagonal matrix\footnote{This structural assumption is imposed in order to use the kernel equations solved in \cite{Deutscher2018tac} in the context of a state feedback design. More general $Q_0$ could be handled by solving new types of kernel equations.}, with arbitrary $Q_1\in\Rset^{n\times n}$.
The dimensions of $F$ and $B$ are implied by those of $\xi(t)$ and $x(z,t)$.
The initial conditions of the system are $\xi(0)=\xi_0\in\Rset^{n_0}$ and $x(z,0)=x_0(z)\in\Rset^n$, $z\in[0,1]$.

Based on the boundary measurement $y(t)\in\Rset^n$ in \eqref{eq:sys-meas}, the objective is to obtain estimates $\hat\xi(t)\in\Rset^{n_0}$ and $\hat x(z,t)\in\Rset^n$ of the states $\xi(t)$ and $x(z,t)$ such that the errors $\varepsilon(t)=\xi(t)-\hat\xi(t)$ and $e(z,t)=x(z,t)-\hat x(z,t)$ go to zero as $t\to\infty$.
For that, the following assumption is imposed.

\begin{assumption}
	\label{ass:detect}
	The pair $(F,C)$ is detectable.
\end{assumption}

As \eqref{eq:sys-rb0} takes the role of an output equation w.r.t.\ the ODE, this necessary condition ensures that at least the unstable part of \eqref{eq:sys-ode} affects the PDE subsystem and can, thus, be detected through the measurement $y(t)$.
Specifically, $C\xi(t)$ acts as a virtual measurement that is subsequently reconstructed through the PDE boundary relation.
This observation motivates the successive design.
Note that the PDE subsystem itself is observable due to its structure and the full boundary measurement at $z=1$.

Consider the observer
\begin{subequations}
	\label{eq:obs}
	\begin{align}
		\dot{\hat\xi}(t) &= F\hat\xi(t) + B\hat x(0,t) + v_\xi(t) \\
		\partial_z \hat x(0,t) &= Q_0\hat x(0,t) + C\hat\xi(t) + v_0(t) \\
		\partial_t \hat x(z,t) &= \Lambda(z)\partial_z^2 \hat x(z,t) + A(z)\hat x(z,t) + v(z,t) \\
		\partial_z \hat x(1,t) &= Q_1 \hat x(1,t) + v_1(t)
	\end{align}
\end{subequations}
with initial conditions $\hat\xi(0)=\hat\xi_0\in\Rset^{n_0}$ and $\hat x(z,0)=\hat x_0(z)\in\Rset^n$, $z\in[0,1]$, for the state estimates $\hat\xi(t)$ and $\hat x(z,t)$.
It is a copy of \eqref{eq:sys} with injection functions $v_\xi(t)\in\Rset^{n_0}$ and $v_0(t),v_1(t),v(z,t)\in\Rset^n$ that may only depend on known quantities, i.e., the measurement \eqref{eq:sys-meas} and the observer state.
This general interpretation of an observer dynamics dates back to the seminal work \cite{Luenberger1966tac} and has shown to be advantageous both for finite-dimensional (nonlinear) systems as well as those of infinite dimension, e.g., in the context of dynamic observers (see \cite{Gehring2025letters}).
When it comes to the backstepping design of observers for distributed-parameter systems, typically, the injection functions are directly chosen to be proportional to the measurement error (see, e.g., \cite{Smyshlyaev2005scl}).
Although this will also be the result here, it is not done at the beginning to better highlight the general design strategy.

In what follows, the injection functions are designed such that the observer error system
\begin{subequations}
	\label{eq:obserror}
	\begin{align}
		\label{eq:obserror-ode}
		\dot\varepsilon(t) &= F\varepsilon(t) + Be(0,t) - v_\xi(t) \\
		\label{eq:obserror-bc0}
		\partial_z e(0,t) &= Q_0e(0,t) + C\varepsilon(t) - v_0(t) \\
		\label{eq:obserror-pde}
		\partial_t e(z,t) &= \Lambda(z)\partial_z^2 e(z,t) + A(z)e(z,t) - v(z,t) \\
		\label{eq:obserror-bc1}
		\partial_z e(1,t) &= Q_1 e(1,t) - v_1(t)
	\end{align}
\end{subequations}
resulting from \eqref{eq:sys} and \eqref{eq:obs} is exponentially stable in an appropriate norm.
For that, a multi-step design makes use of the \emph{strict feedforward form} of the system \eqref{eq:sys} w.r.t.\ the measurement $y(t)$ at $z=1$, a structure that is inherited by the error dynamics \eqref{eq:obserror}.
Roughly speaking, this form implies that $\dot\xi(t)$ (or equivalently $\dot\varepsilon(t)$) may depend on both the ODE and the PDE state, while a change of $x(z,t)$ (or $e(z,t)$) for some specific $z=z_0\in[0,1]$ may only depend on the PDE state itself for $z\in[z_0,1]$.
The strict feedforward form is also apparent when applying a spatial discretization to the PDE subsystem (similar to \cite{Balogh2002ejc}), especially in the case $n=1$, and is highlighted by the ordering of equations in \eqref{eq:sys} and \eqref{eq:obserror}.

\section{Successive backstepping observer design}
\label{sec:obs}

The main idea of the two-step observer design is to exploit the strict feedforward form of the observer error dynamics \eqref{eq:obserror} (and of the system \eqref{eq:sys}) and to construct an observer based on successive transformations.
The first step aims at stabilizing the ODE \eqref{eq:obserror-ode}, which is the subsystem most distal from the measurement at $z=1$, by making use of a virtual measurement.
Together with a standard Volterra integral transformation in the second step, the two successive state transformations essentially move the impact of any destabilizing term to where they can easily be compensated by choice of the injection functions.
In the end, the observer error dynamics is a cascade of respective exponentially stable ODE and PDE subsystems.

\subsection{Step 1: Stabilization of the ODE subsystem via a virtual measurement}
\label{sec:step1}

First, only the subsystem most distal from the measurement $y(t)$ is considered, i.e., the ODE \eqref{eq:sys-ode}.
With regard to that subsystem, the boundary condition \eqref{eq:sys-rb0} or specifically
\begin{equation}
	y_\text{virtual}(t) = C\xi(t)
\end{equation}
takes the role of a virtual measurement.
It represents the impact of the ODE on the PDE subsystem.
To stabilize the ODE error subsystem \eqref{eq:obserror-ode} based on $y_\text{virtual}(t)$, one would want to choose the (virtual) injection
\begin{equation}
	\label{eq:injection-virtual}
	v_\xi(t) = L\big(y_\text{virtual}(t)-C\hat\xi(t)\big) = LC\varepsilon(t),
\end{equation}
at least if the impact of the PDE state is disregarded for the moment.
In view of Assumption~\ref{ass:detect}, it is always possible to find a matrix $L\in\Rset^{n_0\times n}$ in \eqref{eq:injection-virtual} such that $F+LC$ is Hurwitz.
However, the virtual measurement $C\xi(t)$ is obviously not known and, thus, the injection \eqref{eq:injection-virtual} cannot be implemented.
Still, in view of boundary condition \eqref{eq:obserror-bc0} and motivated by \eqref{eq:injection-virtual}, a zero term can be added to the ODE \eqref{eq:obserror-ode} to obtain
\begin{equation}
	\label{eq:obserror-bc0-rewrite}
	\begin{aligned}
		\dot\varepsilon(t) &= F\varepsilon(t) + Be(0,t) - v_\xi(t) \\
		&\qquad + L\big(C\varepsilon(t) - \partial_z e(0,t) + Q_0e(0,t) - v_0(t)\big) \\
		&= (F+LC)\varepsilon(t) + (B+LQ_0)e(0,t) \\
		&\qquad - L\partial_z e(0,t) - v_\xi(t) - Lv_0(t),
	\end{aligned}
\end{equation}
with $F+LC$ Hurwitz.
Importantly, this does not change the dynamics.
In order to express $\partial_z e(0,t)$ in \eqref{eq:obserror-bc0-rewrite} in terms of $e(z,t)$ (and the injection functions), PDE \eqref{eq:obserror-pde} is integrated\footnote{Note that substituting $\partial_z e(0,t)$ by the boundary condition \eqref{eq:obserror-bc0} would void the previous step and simply recover \eqref{eq:obserror-ode}.}:%
\begin{multline}
	\tint_0^1 \partial_t e(z,t)\,\dd z = \big(\Lambda(1)Q_1-\Lambda^\prime(1)\big) e(1,t) \\
	+ \Lambda^\prime(0)e(0,t) + \tint_0^1 \big(\Lambda^\pprime(z)+A(z)\big)e(z,t)\,\dd z \\
	- \Lambda(0)\partial_z e(0,t) - \tint_0^1 v(z,t)\,\dd z - \Lambda(1) v_1(t).
\end{multline}
As using this to replace $\partial_z e(0,t)$ in \eqref{eq:obserror-bc0-rewrite} introduces the time derivative $\partial_t e(z,t)$ in the ODE, a transformation
\begin{equation}
	\label{eq:mot-trafo}
	\bar\varepsilon(t) = \varepsilon(t) - \tint_0^1 L\Lambda^{-1}(0)e(z,t)\,\dd z
\end{equation}
is necessary to recover a state representation
\begin{multline}
	\dot{\bar\varepsilon}(t) = (F+LC)\bar\varepsilon(t) + B_0 e(0,t) \\
	+ \tint_0^1 B(z) e(z,t)\,\dd z + B_1 e(1,t) + \nu(t)
\end{multline}
in the classical sense, with appropriate matrices $B_0$, $B_1$ and $B(z)$ as well as
\begin{equation}
	\nu(t) = \bar L\left[\tint_0^1 v(z,t)\,\dd z + \Lambda(1)v_1(t)\right] - v_\xi(t) - Lv_0(t)
\end{equation}
and $\bar L=L\Lambda^{-1}(0)$.
Provided that the PDE error $e(z,t)$ goes to zero, this construction ensures $\lim_{t\to\infty} \bar\varepsilon(t) = 0$ for an appropriate injection function $\nu(t)$.
The transformation \eqref{eq:mot-trafo} can therefore be interpreted as a systematic procedure for replacing the unavailable virtual measurement $C\xi(t)$ by the available boundary measurement without destroying the desired ODE stabilization.

Motivated by \eqref{eq:mot-trafo}, the first step of the observer design makes use of the (inverse) state transformation
\begin{subequations}
	\label{eq:trafo1}
	\begin{align}
		\varepsilon(t) &= \bar\varepsilon(t) + \tint_0^1 N(z)\bar e(z,t)\,\dd z \\
		e(z,t) &= \bar e(z,t)
	\end{align}
\end{subequations}
with $N(z)\in\Rset^{n_0\times n}$ for $z\in[0,1]$.
Its purpose is to remove the influence of the PDE error from the ODE dynamics and to stabilize the latter.
For that, as in \eqref{eq:mot-trafo}, the condition $N(0)=L\Lambda^{-1}(0)$ has to hold, with $L$ such that $F+LC$ is Hurwitz.
The remaining degrees of freedom in the choice of $N(z)$ are used to decouple \eqref{eq:obserror-ode} from the PDE state $e(z,t)=\bar e(z,t)$, $z\in[0,1)$.
This is done so the observer error dynamics can ultimately be written as the cascade of two stable subsystems.
In order for \eqref{eq:trafo1} to map \eqref{eq:obserror} into the desired form
\begin{subequations}
	\label{eq:trafo1-obserror}
	\begin{align}
		\label{eq:trafo1-obserror-ode}
		\dot{\bar\varepsilon}(t) &= (F+LC)\bar\varepsilon(t) + \bar B \bar e(1,t) - \bar v_\xi(t) \\
		\partial_z \bar e(0,t) &= Q_0 \bar e(0,t) + C\big(\bar\varepsilon(t)+\tint_0^1 N(z)\bar e(z,t)\,\dd z\big) \nonumber\\
		&\hspace{4cm} - v_0(t) \\
		\partial_t \bar e(z,t) &= \Lambda(z) \partial_z^2 \bar e(z,t) + A(z)\bar e(z,t) - v(z,t) \\
		\partial_z \bar e(1,t) &= Q_1\bar  e(1,t) - v_1(t)
	\end{align}
\end{subequations}
with
\begin{align}
	\label{eq:def-barB}
	\bar B &= \dd_z(N\Lambda)(1)-(N\Lambda)(1)Q_1 \\
	\label{eq:def-barvxi}
	\bar v_\xi(t) &= v_\xi(t) - \tint_0^1 N(z)v(z,t)\,\dd z + Lv_0(t) \\
	&\hspace{4cm} - (N\Lambda)(1)v_1(t), \nonumber
\end{align}
the matrix $N(z)$ has to satisfy the initial value problem
\begin{subequations}
	\label{eq:trafo1-ivp}
	\begin{align}
		\dd_z^2(N\Lambda)(z) &= (F+LC)N(z) - N(z)A(z) \\
		(N\Lambda)(0) &= L \\
		\dd_z(N\Lambda)(0) &= B + LQ_0.
	\end{align}
\end{subequations}
It is straightforward to show that \eqref{eq:trafo1-ivp} admits a unique $C^2([0,1])$ solution $N(z)$ (e.g., \cite{Kailath1980}).

It is apparent from \eqref{eq:trafo1-obserror-ode} that choosing
\begin{equation}
	\label{eq:obsgain-1}
	\bar v_\xi(t) = \bar B \bar e(1,t)
\end{equation}
already ensures exponential convergence of the ODE error $\bar\varepsilon(t)$ to zero.
Moreover, the stabilization of the ODE subsystem and its decoupling from the PDE error $\bar e(z,t)$ means that stabilizing the PDE subsystem in the next design step is sufficient to guarantee overall stability.

\subsection{Step 2: Stabilization of the PDE subsystem}

At this point, the first design step has achieved the essential decomposition.
The ODE error evolves as an exponentially stable subsystem driven only by the PDE boundary error \( \bar e(1,t)\) (and the injection function $\bar v_\xi(t)$). Consequently, the remaining task is purely a PDE stabilization problem. This is the reason why a standard Volterra integral transformation%
\begin{subequations}
	\label{eq:trafo2}
	\begin{align}
		\bar\varepsilon(t) &= \tilde\varepsilon(t) \\
		\bar e(z,t) &= \tilde e(z,t) + \tint_z^1 K(z,\zeta)\tilde e(\zeta,t)\,\dd\zeta \eqqcolon \mathcal T^{-1}[\tilde e(t)](z)
	\end{align}
\end{subequations}
can be employed in the second step, where the kernel $K(z,\zeta)\in\Rset^{n\times n}$ is defined on the triangular domain $\mathcal D=\{(z,\zeta)\in[0,1]^2|z\le\zeta\}$.
In order for \eqref{eq:trafo2} to map the observer error dynamics \eqref{eq:trafo1-obserror} into the form
\begin{subequations}
	\label{eq:trafo2-sys}
	\begin{align}
		\dot{\tilde\varepsilon}(t) &= (F+LC)\tilde\varepsilon(t) + \bar B \tilde e(1,t) - \bar v_\xi(t) \\
		\partial_z \tilde e(0,t) &= C\tilde\varepsilon(t) + \tint_0^1 A_0(z)\tilde e(z,t)\,\dd z - v_0(t) \\
		\partial_t \tilde e(z,t) &= \Lambda(z)\partial_z^2 \tilde e(z,t) - \mu \tilde e(z,t) + \!H(z)\tilde e(1,t) - \tilde v(z,t) \\
		\partial_z \tilde e(1,t) &= \big(Q_1+K(1,1)\big) \tilde e(1,t) - v_1(t),
	\end{align}
\end{subequations}
with a strictly upper triangular matrix $A_0(z)\in\Rset^{n\times n}$, a design parameter $\mu$ and (implicit) definitions
\begin{align}
	\mathcal T^{-1}[H](z) &= K(z,1)\Lambda^\prime(1) + \partial_\zeta K(z,1)\Lambda(1) \\
	&\hspace{1cm }- K(z,1)\Lambda(1)\big(Q_1+K(1,1)\big) \nonumber \\
	\label{eq:def-tildev}
	\mathcal T^{-1}[\tilde v(t)](z) &= v(z,t) - K(z,1)\Lambda(1)v_1(t),
\end{align}
the kernel $K(z,\zeta)$ has to satisfy
\begin{subequations}
	\label{eq:trafo2-kerneleqs}
	\begin{align}
		& \Lambda(z)\partial_z^2 K(z,\zeta) - \dd_\zeta^2\big(K(z,\zeta)\Lambda(\zeta)\big) \\
		& \hspace{1cm} = -(A(z)+\mu I)K(z,\zeta) \nonumber \\
		& K(z,z)\Lambda(z) - \Lambda(z)K(z,z) = 0 \\
		& K(z,z)\Lambda^\prime(z) + \partial_\zeta K(z,z)\Lambda(z) + \Lambda(z)\tfrac{\dd}{\dd z}K(z,z) \\
		& \hspace{1cm} + \Lambda(z)\partial_z K(z,z) = A(z) + \mu I \nonumber \\
		\label{eq:trafo2-kerneleqs-0}
		& \partial_z K(0,\zeta) - Q_0 K(0,\zeta) \\
		& \hspace{1cm} = -A_0(\zeta) + CN(\zeta) + \tint_0^\zeta CN(\bar\zeta)K(\bar\zeta,\zeta)\,\dd\bar\zeta \nonumber \\
		& K(0,0) = -Q_0,
	\end{align}
\end{subequations}
where $I$ denotes the identity matrix.
These kernel equations can be traced back\footnote{The explicit calculations become apparent in the context of Section~\ref{sec:dual-step2}.} to those in \cite{Deutscher2018tac}, thus verifying the existence of a unique $C^2(\mathcal D)$ solution of \eqref{eq:trafo2-kerneleqs}.
This solution also defines the non-zero elements of $A_0(z)$ by means of \eqref{eq:trafo2-kerneleqs-0} (see again \cite{Deutscher2018tac}).

The structure of the transformed system \eqref{eq:trafo2-sys} stands out in that all destabilizing elements of the coupled PDE-ODE system are contained in terms proportional to the measurement error $\tilde e(1,t)$.
This greatly facilitates the choice of stabilizing injection functions because a simple compensation is sufficient.
As \eqref{eq:obsgain-1} already ensures exponential stability of the ODE subsystem, setting
\begin{subequations}
	\label{eq:obsgain-2}
	\begin{align}
		v_0(t) &= 0 \\
		v_1(t) &= \big(Q_1+K(1,1)\big) \tilde e(1,t) \\
		\tilde v(z) &= H(z) \tilde e(1,t)
	\end{align}
\end{subequations}
guarantees convergence of the PDE state $\tilde e(z,t)$ to zero, provided the design parameter $\mu$ is chosen appropriately.
Thus, the observer error dynamics
\begin{subequations}
	\label{eq:trafo2-obs-error}
	\begin{align}
		\dot{\tilde\varepsilon}(t) &= (F+LC)\tilde\varepsilon(t) \\
		\partial_z \tilde e(0,t) &= C\tilde\varepsilon(t) + \tint_0^1 A_0(z)\tilde e(z,t)\,\dd z \\
		\partial_t \tilde e(z,t) &= \Lambda(z)\partial_z^2 \tilde e(z,t) - \mu \tilde e(z,t) \\
		\partial_z \tilde e(1,t) &= 0
	\end{align}
\end{subequations}
obtained from \eqref{eq:trafo2-sys} with \eqref{eq:obsgain-1} and \eqref{eq:obsgain-2} is exponentially stable in an appropriate norm, which is detailed next.


\subsection{State observer}

Let $L$ be such that $F+LC$ is Hurwitz, which is always possible by Assumption~\ref{ass:detect}, and choose $\mu>0$.
Then, the observer error dynamics \eqref{eq:trafo2-obs-error} is exponentially stable in the norm $||\cdot||=(||\cdot||_{\Rset^{n_0}}^2+||\cdot||_{L_2^n}^2)^\frac{1}{2}$ for all initial conditions $\varepsilon(0)$ and $e(z,0)$, $z\in[0,1]$ (see, e.g., \cite{Deutscher2021TAC}).
This implies the same for \eqref{eq:obserror} by the obvious invertibility of the transformations \eqref{eq:trafo1} and \eqref{eq:trafo2}.
The injection functions
\begin{subequations}
	\label{eq:obsgain-all}
	\begin{align}
		v_\xi(t) &= \Big[\dd_z(N\Lambda)(1) + (N\Lambda)(1)K(1,1) \\
		&\hspace{0cm} + \! \tint_0^1 \! N(z) \big(K(z,1)\Lambda^\prime(1)+\partial_\zeta K(z,1)\Lambda(1)\big)\dd z\Big] \! e(1,t) \nonumber \\
		v_0(t) &= 0, \qquad v_1(t) = \big(Q_1+K(1,1)\big) e(1,t) \\
		v(z,t) &= \big(K(z,1)\Lambda^\prime(1)+\partial_\zeta K(z,1)\Lambda(1)\big) e(1,t)
	\end{align}
\end{subequations}
for the observer \eqref{eq:obs} result from the choice in \eqref{eq:obsgain-1} and \eqref{eq:obsgain-2} by using \eqref{eq:def-barB}, \eqref{eq:def-barvxi} and \eqref{eq:def-tildev} in light of $\tilde e(1,t)=\bar e(1,t)=e(1,t)$.
Note that these injection functions only depend on the measurement error, i.e., on known quantities.

\section{Duality to controller design}
\label{sec:dual}

The duality developed in this section serves two purposes.
First, it establishes that the proposed observer design is the dual counterpart of a multi-step state feedback design.
Second, and more importantly, it explains why the observer can and should be constructed successively.
The strict feedforward structure of the primal system becomes a strict feedback structure in the dual system, so that the observer's first step corresponds to the virtual control step of the state feedback design.
Duality is not just shown between the respective transformations and corresponding kernel equations as well as the systems themselves but also between the error injection functions and the state feedback.
To illustrate the relation, a rough notion of duality in the context of parabolic PDE-ODE systems is presented, without going into mathematical details.
The two-step state feedback design used here constitutes a simplified version of the approach in \cite{Gehring2021MTNS}.

\subsection{Duality and dual system}

From a mathematical point of view, a system dynamics $(\dot X(t),Y(t)) = \Sigma(X(t),U(t))$ can be understood as a map $\Sigma:\mathcal X\times \mathcal U\to\mathcal X\times \mathcal Y$ from the current state $X(t)\in\mathcal X$ and input $U(t)\in\mathcal U$ to the state derivative $\dot X(t)\in\mathcal X$ and output $Y(t)\in\mathcal Y$.
Denoting by $\langle\cdot,\cdot\rangle$ the inner product associated with the corresponding function spaces, the following definition of a dual system is based on the notation of adjoint operators (see, e.g., \cite{Curtain2020}).

\begin{definition}
	\label{def:duality}
	Let $X$, $U$ and $Y$ be the state, input and output of a primal system $(\dot X,Y) = \Sigma(X,U)$.
	Then, the associated dual system $(\dot X^\ast,Y^\ast) = \Sigma^\ast(X^\ast,U^\ast)$ with state $X^\ast$, input $U^\ast$ and output $Y^\ast$ satisfies
	\begin{equation}
		\label{eq:duality}
		\langle \Sigma(X,U), (X^\ast, U^\ast)\rangle \stackrel{!}{=} \langle (X,U), \Sigma^\ast(X^\ast, U^\ast)\rangle
	\end{equation}
\end{definition}

Here, the natural choice $\mathcal X = \Rset^{n_0}\times (L^2([0,1]))^n$ of a state space for \eqref{eq:sys} implies the (weighted) inner product
\begin{multline}
	\langle X_1(t),X_2(t)\rangle_{\mathcal X} = \tint_0^1 x_1^T(z,t)\Lambda^{-1}(z)x_2(z,t)\,\dd z \\
	+ \xi_1^T(t)\xi_2(t)
\end{multline}
for two elements $X_1(t),X_2(t)\in\mathcal X$.
In what follows, the weight $\Lambda^{-1}(z)$ for $L^2([0,1])$ functions is introduced for simplification only, while the standard scalar product is adopted for all real-valued vectors.
By applying \eqref{eq:duality} to the dynamics \eqref{eq:sys} with $X(t)=(\xi(t),x(\cdot,t))\in\mathcal X$ and $Y(t)=y(t)\in\Rset^n$, the dual system
\begin{subequations}
	\label{eq:dual-sys}
	\begin{align}
		\dot\xi^\ast(t) &= F^\ast\xi^\ast(t) + B^\ast x^\ast(0,t) \\
		\partial_z x^\ast(0,t) &= Q_0^\ast x^\ast(0,t) + C^\ast\xi^\ast(t) \\
		\partial_t x^\ast(z,t) &= \Lambda(z)\partial_z^2 x^\ast(z,t) + A^\ast(z)x^\ast(z,t) \\
		\partial_z x^\ast(1,t) &= Q_1^\ast x^\ast(1,t) + u^\ast(t)
	\end{align}
\end{subequations}
is obtained, where
\begin{equation}
	\label{eq:def-Xast}
	X^\ast(t) = (\xi^\ast(t),x^\ast(\cdot,t))\in\mathcal X,
\end{equation}
and $U^\ast(t)=u^\ast(t)\in\Rset^n$ as well as
\begin{align}
	\label{eq:dual-defs}
	F^\ast &= F^T, & B^\ast &= -C^T, & C^\ast &= -B^T \\
	Q_0^\ast &= Q_0^T, & Q_1^\ast &= Q_1^T, & A^\ast(z) &= \Lambda(z)A^T(z)\Lambda^{-1}(z). \nonumber
\end{align}
As \eqref{eq:sys} is a system without input, \eqref{eq:dual-sys} does not have an output.
The strict feedforward form of \eqref{eq:sys} w.r.t.\ the output $y(t)$ results in a strict feedback form of the dual system \eqref{eq:dual-sys} w.r.t.\ the input $u^\ast(t)$.
In particular, detectability of $(F,C)$ (recall Assumption~\ref{ass:detect}) implies stabilizability of $(F^\ast,B^\ast)$.

The stabilization of \eqref{eq:dual-sys} by means of a state feedback for $u^\ast(t)$ is dual to the stabilization of the error dynamics \eqref{eq:obserror} by choice of $v_\xi(t)$, $v_0(t)$, $v_1(t)$ and $v(z,t)$.
In both cases, the goal is to drive the state to zero.
As the feedback design relies on knowledge of the system state, duality of \eqref{eq:obserror} with
\begin{equation}
	\label{eq:def-X}
	X(t) = (\varepsilon(t),e(\cdot,t))\in\mathcal X,
\end{equation}
and \eqref{eq:dual-sys} with $X^\ast(t)$ in \eqref{eq:def-Xast} can be shown by defining the known quantities as measurements:
\begin{multline}
	\label{eq:def-Yast}
	Y^\ast(t) = (-\xi^\ast(t),x^\ast(0,t),-x^\ast(1,t),-x^\ast(\cdot,t)) \\
	\in\Rset^{n_0+2n}\times L^2([0,1]).
\end{multline}
For that, set 
$U(t)=(v_\xi(t),v_0(t),v_1(t),v(\cdot,t))\in\Rset^{n_0+2n}\times L^2([0,1])$ as well as $Y(t)=e(1,t)\in\Rset^n$ based on \eqref{eq:obserror} and apply \eqref{eq:duality} in view of \eqref{eq:dual-defs}.
Note that the signs in \eqref{eq:def-Yast} are introduced only to show duality.

\begin{remark}
	Duality of \eqref{eq:obserror} and \eqref{eq:dual-sys} can also be shown by defining the injection functions to be proportional to the measurement error $e(1,t)$ on the one hand and the state feedback as a linear combination of $\xi^\ast(t)$, $x^\ast(0,t)$, $x^\ast(1,t)$ and an integral over $x^\ast(z,t)$ on the other hand, and applying \eqref{eq:duality} in view of both systems having neither input nor output.
\end{remark}

\subsection{1st design step}
\label{sec:dual-step1}

Following the design in \cite{Gehring2021MTNS}, the boundary value $x^\ast(0,t)$ takes the role of a (virtual) control input w.r.t.\ the ODE subsystem in view of the strict feedback form of \eqref{eq:dual-sys}.
Thus, one would choose the virtual state feedback $x^\ast(0,t) = K^\ast\xi^\ast(t)$ such that $F^\ast+B^\ast K^\ast$ is Hurwitz.
This is obviously not possible and instead implies a state transformation
\begin{subequations}
	\label{eq:dual-trafo1}
	\begin{align}
		\bar\xi(t) &= \xi(t) \\
		\bar x^\ast(z,t) &= x^\ast(z,t) - N^\ast(z)\xi^\ast(t)
	\end{align}
\end{subequations}
with the ODE state left unchanged and $N^\ast(z)\in\Rset^{n\times n_0}$, $z\in[0,1]$, as the solution of the initial value problem
\begin{subequations}
	\label{eq:dual-trafo1-ivp}
	\begin{align}
		\Lambda(z)\dd_z^2 N^\ast(z) &= N^\ast(z)(F^\ast+B^\ast K^\ast) - A^\ast(z)N^\ast(z) \\
		N^\ast(0) &= K^\ast \\
		\dd_z N^\ast(0) &= C^\ast + Q_0^\ast K^\ast,
	\end{align}
\end{subequations}
in order to map \eqref{eq:dual-sys} into the form
\begin{subequations}
	\label{eq:dual-trafo1-sys}
	\begin{align}
		\dot{\bar\xi}^\ast(t) &= (F^\ast+B^\ast K^\ast)\bar\xi^\ast(t) + B^\ast \bar x^\ast(0,t) \\
		\partial_z \bar x^\ast(0,t) &= Q_0^\ast \bar x^\ast(0,t) \\
		\partial_t \bar x^\ast(z,t) &= \Lambda(z) \partial_z^2 \bar x^\ast(z,t) + A^\ast(z) \bar x^\ast(z,t) \\
		&\hspace{3cm} - N^\ast(z) B^\ast \bar x^\ast(0,t) \nonumber \\
		\label{eq:dual-trafo1-sys-bc1}
		\partial_z \bar x^\ast(1,t) &= Q_1^\ast \bar x^\ast(1,t) + \bar C^\ast \bar\xi^\ast(t) + u^\ast(t)
	\end{align}
\end{subequations}
with $\bar C^\ast = Q_1^\ast N^\ast(1)- \dd_z N^\ast(1)$, where the ODE is exponentially stable (if $\bar x^\ast(0,t)=0$).

This design step is dual to that of the observer design in Section~\ref{sec:step1}, which uses a virtual measurement.
To verify that, define by $X(t) = \mathcal T_1^{-1}[\bar X](t)$ the transformation in \eqref{eq:trafo1} and by $\bar X^\ast(t)=\mathcal T_1^\ast[X^\ast](t)$ that in \eqref{eq:dual-trafo1}, with $\mathcal X(t)$ as in \eqref{eq:def-X}, $\mathcal X^\ast(t)$ as in \eqref{eq:def-Xast} as well as
\begin{subequations}
	\label{eq:def-barX}
	\begin{align}
		\bar{\mathcal X}(t) &= (\bar\varepsilon(t),\bar e(\cdot,t))\in\mathcal X \\
		\bar{\mathcal X}^\ast(t) &= (\bar\xi^\ast(t),\bar x^\ast(\cdot,t))\in\mathcal X.
	\end{align}
\end{subequations}
Then, similar to \eqref{eq:duality}, both transformations are dual in the sense of
\begin{equation}
	\label{eq:duality-trafo}
	\langle \mathcal T^{-1}[\bar X], X^\ast \rangle \stackrel{!}{=} \langle \bar X, \mathcal T^\ast[X^\ast] \rangle
\end{equation}
if
\begin{equation}
	\label{eq:dual-def-Nast}
	N^\ast(z) = -\Lambda(z)N^T(z).
\end{equation}
In fact, substituting the latter relation in \eqref{eq:dual-trafo1-ivp} and setting $K^\ast=-L^T$, the initial value problem \eqref{eq:trafo1-ivp} for $N(z)$ reemerges.
It is therefore unsurprising that duality of the transformed systems \eqref{eq:trafo1-obserror} and \eqref{eq:dual-trafo1-sys} can be shown using \eqref{eq:duality} with $U(t)=(\bar v_\xi(t),v_0(t),v_1(t),v(\cdot,t))\in\Rset^{n_0+2n}\times L^2([0,1])$, $Y(t)=e(1,t)\in\Rset^n$, $U^\ast(t)=u^\ast(t)\in\Rset^n$ and \eqref{eq:def-Yast}.
For that, \eqref{eq:dual-defs}, \eqref{eq:dual-def-Nast}, $K^\ast=-L^T$ and $\bar C^\ast=\bar B^T$ are taken into account.

\begin{remark}
	In the context of adjoint operators, it should be $(\mathcal T^\ast)^{-1}$ instead of $\mathcal T^\ast$ in \eqref{eq:duality-trafo}.
	However, as the transformations \eqref{eq:trafo1} and \eqref{eq:trafo2} in the observer design are inverse maps and the dual ones in \eqref{eq:dual-trafo1} and \eqref{eq:dual-trafo2} for the control design are not, the definition \eqref{eq:duality-trafo} simplifies the interpretation of the relation between the primal and the dual results.
\end{remark}

\subsection{2nd design step}
\label{sec:dual-step2}

The Volterra integral transformation
\begin{subequations}
	\label{eq:dual-trafo2}
	\begin{align}
		\tilde\xi^\ast(t) &= \bar\xi^\ast(t) \\
		\tilde x^\ast(z,t) &= \bar x^\ast(z,t) - \tint_0^z K^\ast(z,\zeta)\bar x^\ast(\zeta,t)\,\dd\zeta
	\end{align}
\end{subequations}
with $K^\ast(z,\zeta)\in\Rset^{n\times n}$ on $\mathcal D^\ast=\{(z,\zeta)\in[0,1]^2|\zeta\le z\}$ satisfying the kernel equations
\begin{subequations}
	\label{eq:dual-trafo2-kerneleqs}
	\begin{align}
		& \Lambda(z)\partial_z^2 K^\ast(z,\zeta) - \dd_\zeta^2\big(K^\ast(z,\zeta)\Lambda(\zeta)\big) \\
		&\hspace{1cm} = K^\ast(z,\zeta)\big(A^\ast(\zeta)+\mu I\big) \nonumber \\
		& K^\ast(z,z)\Lambda(z) - \Lambda(z)K^\ast(z,z) = 0 \\
		& K^\ast(z,z)\Lambda^\prime(z) + \partial_\zeta K^\ast(z,z)\Lambda(z) + \Lambda(z)\tfrac{\dd}{\dd z}K^\ast(z,z) \! \\
		&\hspace{1cm} \Lambda(z)\partial_z K^\ast(z,z) = -\big(A^\ast(z)+\mu I\big) \nonumber \\
		& \partial_\zeta K^\ast(z,0)\Lambda(0) + K^\ast(z,0)\big(\Lambda^\prime(0)-\Lambda(0)Q_0^\ast\big) \\
		&\hspace{1cm} = -A_0^\ast(z) - N^\ast(z)B^\ast + \tint_0^z K^\ast(z,\zeta)N^\ast(\zeta)B^\ast\,\dd\zeta \nonumber \\
		& K^\ast(0,0) = Q_0^\ast
	\end{align}
\end{subequations}
maps \eqref{eq:dual-trafo1-sys} into the form
\begin{subequations}
	\label{eq:dual-trafo2-sys}
	\begin{align}
		\dot{\tilde\xi}^\ast(t) &= (F^\ast+B^\ast K^\ast)\tilde\xi^\ast(t) + B^\ast\tilde x^\ast(0,t) \\
		\partial_z \tilde x^\ast(0,t) &= 0 \\
		\partial_t \tilde x^\ast(z,t) &= \Lambda(z) \partial_z^2 \tilde x^\ast(z,t) - \mu \tilde x^\ast(z,t) + A_0^\ast(z)\tilde x^\ast(0,t) \\
		\label{eq:dual-trafo2-sys-bc1}
		\partial_z \tilde x^\ast(1,t) &= \big(Q_1^\ast-K^\ast(1,1)\big) \tilde x^\ast(1,t) + \bar C^\ast \tilde\xi^\ast(1) + u^\ast(t) \nonumber \\
		&\hspace{0cm} + \tint_0^1 \big(H^\ast(z) + \tint_z^1 H^\ast(\zeta)K_\text{I}^\ast(\zeta,z)\,\dd\zeta\big) \tilde x^\ast(z,t)\, \dd z
	\end{align}
\end{subequations}
with a strictly lower triangular matrix $A_0^\ast(z)\in\Rset^{n\times n}$ and $H^\ast(z)=\big(Q_1-K^\ast(1,1)\big)K^\ast(1,z)-\partial_z K^\ast(1,z)$, where $K_\text{I}(z,\zeta)\in\Rset^{n\times n}$ is the kernel of the inverse transformation%
\begin{equation}
	\bar x^\ast(z,t) = \tilde x^\ast(z,t) + \tint_0^z K_\text{I}^\ast(z,\zeta)\tilde x^\ast(\zeta,t)\,\dd\zeta
\end{equation}
of \eqref{eq:dual-trafo2}.
Stabilization is achieved by setting the design parameter $\mu$ sufficiently large and $u^\ast(t)$ such that $\partial_z\tilde x^\ast(1,t)=0$.
The invertibility of the transformations \eqref{eq:dual-trafo1} and \eqref{eq:dual-trafo2} then implies that $\xi(t)$ and $x(z,t)$ go to zero if the feedback
\begin{multline}
	\label{eq:dual-u}
	u^\ast(t) = - \big(Q_1^\ast-K^\ast(1,1)\big) x^\ast(1,t) - \big(K^\ast(1,1)N^\ast(1) \\
	- \dd_z N^\ast(1) + \tint_0^1 \partial_z K^\ast (1,z) N^\ast(z)\,\dd z\big) \xi^\ast(t) \\
	+ \tint_0^1 \partial_z K^\ast(1,z)x^\ast(z,t)\,\dd z 
\end{multline}
of the original state is applied to \eqref{eq:dual-sys}.

The transformation $\bar X(t)=\mathcal T_2^{-1}[\tilde X](t)$ associated with \eqref{eq:trafo2} can be shown to be dual in the sense of \eqref{eq:duality-trafo} to $\tilde X^\ast(t) = \mathcal T_2^\ast[\bar X^\ast](t)$, which follows from \eqref{eq:dual-trafo2}, by introducing
\begin{subequations}
	\begin{align}
		\tilde{\mathcal X}(t) &= (\tilde\varepsilon(t),\tilde e(\cdot,t))\in\mathcal X \\
		\tilde{\mathcal X}^\ast(t) &= (\tilde\xi^\ast(t),\tilde x^\ast(\cdot,t))\in\mathcal X
	\end{align}
\end{subequations}
in analogy to \eqref{eq:def-barX} and setting
\begin{equation}
	K^\ast(z,\zeta) = -\Lambda(z)K^T(\zeta,z)\Lambda^{-1}(\zeta).
\end{equation}
Substituting that as well as \eqref{eq:dual-def-Nast} and \eqref{eq:dual-defs} into \eqref{eq:dual-trafo2-kerneleqs}, following some tedious calculations, the kernel equations for the controller design turn out to be identical to those in \eqref{eq:trafo2-kerneleqs} for the observer design if $A_0^\ast(z) = -\Lambda(z)A_0^T(z)$ is set.
This proves that the existence of a unique $C^2(\mathcal D)$ solution $K(z,\zeta)$ for \eqref{eq:trafo2-kerneleqs} is directly implied by the result in \cite{Deutscher2018tac}, where kernel equations \eqref{eq:dual-trafo2-kerneleqs} are considered.
Similar observer-controller duality at the level of kernel PDEs is reported, e.g., in \cite{Smyshlyaev2005scl}.

Duality also holds true for the observer error dynamics \eqref{eq:trafo2-obs-error} and the closed loop \eqref{eq:dual-trafo2-sys} with $\partial_z\tilde x^\ast(1,t)=0$ instead of \eqref{eq:dual-trafo2-sys-bc1}.
To verify that, simply apply \eqref{eq:duality} in view of the relation between the matrices in the primal and dual system.
Finally, close inspection of the state feedback \eqref{eq:dual-u} and the injection functions in \eqref{eq:obsgain-all} reveals duality of the two.

\section{Concluding remarks}

The main result of this work is a design principle for observers of strictly feedforward systems, discussed here for parabolic PDE-ODE systems.
Duality with multi-step state feedback in \cite{Gehring2021MTNS} shows that the observer transformations are not ad hoc constructions but the dual counterparts of the transformations underlying the controller design.
The proposed design strategy can directly be transferred to other strictly feedforward systems, including a more general form%
\begin{subequations}
	\label{eq:sys-general}
	\begin{align}
		\dot\xi(t) &= F\xi(t) + Bx(0,t) + \tint_0^1 G_1(z)x(z,t)\,\dd z \nonumber \\
		&\hspace{.5cm} + G_2x(1,t) + B_1u(t) \\
		\partial_z x(0,t) &= Q_0 x(0,t) + C\xi(t) + \tint_0^1 G_3(z)x(z,t)\,\dd z \nonumber \\
		&\hspace{.5cm} + G_4x(1,t) + B_2u(t) \\
		\partial_t x(z,t) &= \Lambda(z)\partial_z^2 x(z,t) + A(z)x(z,t) + B_3(z)u(t) \nonumber \\
		&\hspace{.5cm} + \tint_z^1 G_5(z,\zeta)x(\zeta,t)\,\dd\zeta + G_6(z)x(1,t) \\
		\partial_z x(1,t) &= Q_1 x(1,t) + u(t).
	\end{align}
\end{subequations}
of \eqref{eq:sys}, potentially with other types of boundary conditions, parabolic ODE-PDE-ODE systems, and hyperbolic distributed-parameter systems with finite-dimensional dynamics at the boundaries.
Importantly, a strict feedforward form w.r.t.\ the boundary measurement is always required.
For instance, \eqref{eq:sys-ode}--\eqref{eq:sys-rb1} is not in strict feedforward form w.r.t.\ $x(0,t)$ or a measurement of the ODE state $\xi(t)$.

\addtolength{\textheight}{-12cm}   


\section*{Acknowledgment}
This research was funded in part by the Austrian Science Fund (FWF) [I 6519-N]. The author also thanks Dr.\ Abdurrahman Irscheid of Saarland University, Germany, for many valuable discussions.

\bibliographystyle{plain}
\bibliography{/home/ngehring/Documents/90_Literatur/mybib}

\end{document}